\documentstyle[]{article}
\begin{document}
  \title{VECTOR FIELDS AND DIFFERENTIAL OPERATORS:\\ NONCOMMUTATIVE CASE.}
\author{
Andrzej Borowiec 
\\Institute of Theoretical Physics,
University of Wroc{\l}aw\\ 
Pl. Maxa Borna 9, 50-204 Wroc{\l}aw, Poland\\ 
e-mail: borow@ift.uni.wroc.pl}
\maketitle
\def\ba{\begin{array}}
\def\ea{\end{array}}
\def\lra{\longrightarrow}
\def\ra{\rightarrow}
\def\ld{\ldots}
\def\cd{\cdot}
\def\be{\begin{equation}}
\def\ee{\end{equation}}
\def\bm{\begin{em}}
\def\en{\end{em}}
\def\bk{I \kern-.25em k} 
\def\N{I \kern-.25em N} 
\def\R{I \kern-.25em R} 
\def\lin{\hbox{lin}\,}
\def\der{\hbox{der}\,}
\def\derk{\der_{\bk}\,}
\def\bimod{\hbox{bimod}\,}
\def\im{\hbox{im}\,}
\def\ker{\hbox{ker}\,}
\def\ann{\hbox{ann}\,}
\def\bot{\,_b\!\otimes}
\def\nd{\hbox{End}\,}
\def\hom{\hbox{Hom}\,}
\def\alg{\hbox{alg}\,}
\def\D{{\cal D}}

\def\id{\hbox{id}}
\def\idA{\id_A}
\def\idV{\id_V}
\def\ot{\otimes}
\def\c{\circ}
\def\bot{\,_b\!\otimes}
\def\l{\lambda}
\def\lc#1{#1^{\lambda}}
\def\p{\partial}
\def\rp#1{#1^{\p}}
\def\lp#1{^{\p}\!#1}
\def\O{\Omega^1_u (A)}
\def\M{\cal M}
\def\x{\cal X}
\def\X{{\cal X}_u(A)}

\def\normalbaselines{\baselineskip20pt \lineskip3pt \lineskiplimit20pt}
\def\mapright#1{\smash{\mathop{\longrightarrow}\limits^{#1}}}
\def\mapleft#1{\smash{\mathop{\longleftarrow}\limits^{#1}}}
\def\mapdown#1{\Big\downarrow\rlap
 {$\vcenter{\hbox{$\scriptstyle#1$}}$}}
\def\mapup#1{\Big\uparrow\rlap
 {$\vcenter{\hbox{$\scriptstyle#1$}}$}}
\def\mapupdown#1{\Big\updownarrow\rlap
 {$\vcenter{\hbox{$\scriptstyle#1$}}$}}

\begin{abstract}
A notion of Cartan pairs as an analogy of vector fields
in the realm of noncommutative geometry has been proposed in \cite{AB}.
In this paper we give an outline of the construction of a noncommutative 
analogy of the algebra of partial differential operators as well as its 
natural (Fock type) representation. We shall also define 
co-universal vector fields and covariant derivatives. 
\end{abstract}
\section{Introduction.}

Let $A$ be an associative (not necessarily commutative) algebra with unit.
A new notion of right (left) {\it Cartan pairs} over algebra $A$ has been 
proposed in our previous paper \cite{AB} as a noncommutative substitute of 
concept of vector fields.  A Cartan pair consist of algebra $A$ and an 
$A$-bimodule $M$ equipped with a suitable right (or left) action of $M$ 
on $A$. As a next step we have  defined
a notion of right (left) dual of an $A$-bimodule  in such
a way that dual object is again an $A$-bimodule.
Finally, it has been shown that the first order differential calculi on 
$A$ and  right (left) $A$-Cartan pairs are dual each other. Therefore, in 
noncommutative settings, the concept of vector fields splits into two
concepts: right and left vector fields.
For modules over commutative algebras  both notions of left and right 
vector fields coincide.
In this case, above dualities restore  classical dualities between one
forms and vector fields which are known from  the classical differential 
geometry on manifolds ($A=C^\infty ({\cal M})$). 

The construction employs the Cartan formula
$$
X(f)\equiv <X,\, df>\equiv \hbox{i}_X df \ . \eqno (1.1)
$$
which expresses both dualities. Our main result was that formula (1.1) 
allows to reconstruct the action of vector fields  on "functions" 
({\it i.e.} elements of the algebra) if we are given the differential.  
Conversely, one can find out the differential 
by means of action. Examples of such actions for  given (noncommutative)
calculi can be found in \cite{BK,DMH,Ma}.

The Leibniz rule for a first-order differential  calculus
$$
d(fg) = (df).\,g + f.\,dg  \eqno (1.2)
$$
remains unchanged as in the classical case:
this one for vector fields (partial derivatives) has been replaced by 
more general axiom:  action of bimodule on the algebra
(see Definition 2.3 below). 

Our formalism is inspired by but different from  the 
{\it Lie--Cartan pairs} approach \cite{KS}.
The last notion (with different names) has been rediscovered 
independently by many authors from time to time (see \cite{Mac}).
In particular, we have no analogue of Lie bracket. Another generalization 
which is based on a derivation property of vector fields (so called 
{\it Lie pairs}) has been introduced in \cite{FGV}.

Our aim here is to push further a {\it contravariant} (vector fields)
formalism in the noncommutative differential geometry. The following 
topics will be discussed:
co-universal vector fields,
algebra of differential operators and
covariant differentiations.

For the sake of brevity we shall develop a "right handed" version of the
theory. It will become soon clear that the similar 
construction can be performed for left Cartan pairs (c.f \cite{AB}).

\section{Preliminaries and notation.}

Throughout this paper $\bk$ denotes some fixed unital and commutative 
ring. (For simplicity one can limit ourselves to the field of real or 
complex numbers.) Algebras are unital associative $\bk$-algebras and 
homomorphisms are assumed to be unital. 
All objects considered here  are $\bk$-modules,
all maps are assumed to be $\bk$-linear maps. 
The tensor product $\ot$ unless otherwise mentioned 
means $\ot_{\bk}$\,.
Let $M$ be an $(A, A)$--bimodule ($A$--bimodule in short). We shall
denote by dot "." the both: left and right multiplication by elements
from $A$. For example, by bimodule axioms, one has $(f.x).g = f.(x.g) = 
f.x.g$ for $f, g\in A$ and $x\in M$. In this Section we give a necessary
definitions and  summarize main results from \cite{AB}.

Let $M$ be an $A$--bimodule and let $\hom_{(-,A)}(M,\, A)$ 
denotes a set of all right $A$-module maps from $M$ into $A$.
It is a $\bk$-space. For $X\in \hom_{(-,A)}(M,\, A)$ denote
as a pairing
$$
<X,\, m>\,\doteq\, X(m) \in A \eqno (2.1)$$
the evaluation of $X$ on the element $m\in M$. 
(see e.g. \cite{Bou} p. 232 in the context of modules over 
commutative algebras.)
Then the following formula
$$
<f.X.g,\, m>\, \doteq\, f<X,\, g.m> \eqno (2.2)$$
where $f,g \in A$,
defines an $A$-bimodule structure on $\hom_{(-,A)}(M,\, A)$ \cite{AB}.
Therefore, the pairing (2.1) defines a bimodule map
$$
<\, ,\,> :\, M^*\ot_A  M \ra A \eqno (2.3)$$

{\bf Definition 2.1}. The set $M^* \doteq \hom_{(-,A)}(M,\, A)$ equipped 
with a bimodule structure indicated above is called a 
{\it right dual} of $M$.  

In a similar way one defines a {\it left dual} 
$\,^* M\,\doteq\, \hom_{(A,-)}(M,\, A)$. These notions generalize
the standard notion of duality: for modules over commutative algebra the 
three notions of duals; right, left and standard dual of module coincide.

{\bf Example 2.2}. Observe that for the algebra itself $A^* = A$ as 
$A$-bimodules. An identification is made by $f\mapsto f(1)$. 
In this case $<f,\, g>=fg$. Of course, in the similar way $^*A=A$.

By an {\it action} of $M$ on $A$ we mean a $\bk$--linear mapping 
$\beta\in\hom_{\bk}(M,\, \nd_{\bk}(A))$. We shall also write
$M\ni x\mapsto x^\beta\in\nd_{\bk}(A)$ or
$A\ni f\mapsto x^{\beta}(f)\in A$ to denote the action $\beta$.

{\bf Definition 2.3}. 
An $A$-{\it right Cartan pair} $(N, \rho)$ consists of an $A$--bimodule 
$N$ and its right action $\rho :N\ra \nd_{\bk}(A)$, such that
$$
(f.X)^\rho (g)\ =\ f\, X^\rho (g) \eqno (2.4)$$
and
$$
X^\rho (f\, g)\ =\ X^\rho (f)\,g + (X.f)^\rho (g) \eqno (2.5)$$
for all $X\in N$. Such action will be called a {\it right action}. 

It is easy to see that $X^\rho$ annihilates scalars from $\bk$.

Observe that in the case of modules over commutative algebras 
$X.f=f.X$ \, and the formulae (2.4) - (2.5)  reduce to the 
standard Leibniz rule 
$$X^\rho (f\, g)\ =\ X^\rho (f)\,g + f\, X^\rho (g) \eqno (2.6)$$
satisfied in the classical differential geometry. Therefore, $X$
becomes a derivation of the algebra $A$.
In this sense the definition of Cartan pairs is a generalization of
notion of vector fields. As we will see below the dualization of 
Cartan  pairs leads to a differential calculus on the algebra $A$.

{\bf Definition 2.4}. A (first order) {\it differential calculus} 
$(M\,, d)$ on the algebra  $A$ (or in short an $M$-{\it valued calculus 
on} $A$) consists of an $A$-bimodule $M$ and a 
linear map $d: A\ra M$ satisfying the Leibniz rule (1.2).

The  bimodule $M$ plays the role of a bimodule of one-forms.
Noncommutative differential calculi {\it i.e.} differential calculi on 
noncomutative algebras are basic objects of noncommutative differential 
geometry \cite{Co}. They have been investigated by many authors 
(see {\it e.g.} \cite{BK,DMH,Ma,PW}).

Let now $(M, d)$ be a calculus on an algebra $A$. The differential
$d$ and formula (1.1) defines an action of the right dual $M^*$ 
on $A$. This action 
$$A\ni f\mapsto\  \rp X(f)\ \doteq\  <X,\, d f> \eqno (2.7)$$
will be called a {\it right partial derivatives} along the element $X\in 
M^*$ with respect to the calculus $(M, d)$. It appears that this action
satisfies axioms of right Cartan pairs. Indeed, by (2.2)
$$
\rp{(f.X)} (g)= <f.X,\, dg> = f\,<X,\, dg> = f\rp X(g)$$
and
$$
\rp{(X.f)} (g)= <X.f,\, dg> = <X,\, f.dg> = 
<X,\, d(fg)- df.g> = \rp X(fg) - \rp X(f) g \ . $$
Therefore, to each differential
calculus $(M, d)$ on $A$ we can  associate a unique right 
Cartan pair of right partial derivatives $(M^*, \p)$ of $(M, d)$. The 
converse statement is also true: to each right Cartan pair $(N, \rho)$ one 
can associate a unique differential calculus $(^*\!N, d_\rho)$ where,
$d_\rho:A\ra\ ^*\!N$  is defined by formula (2.8) below. Thus we have

{\bf Theorem 2.5}. {\em
Let $(M, d)$ be a calculus on $A$. Then $M^*$ together with an action (2.7),
via the right partial derivatives, becomes the right Cartan pair 
$(M^*, \p)$ on $A$. 

Conversely, let $(N, \rho)$ be a right Cartan pair on $A$. Then the formula
$$
<X, d_\rho f>\  =\  X^\rho (f) \eqno (2.8)$$
for each $X\in N$, determines $d_\rho f$ as an element of a left dual 
$^*\!N$ of the bimodule $M$. The mapping $d_\rho : A\ra\ ^*\!N$ defines 
the $^*\!N$--valued calculus $(^*\!N, d_\rho)$ on $A$.

Moreover, the module of one forms is spanned by 
differential if and only if  the action has a trivial kernel.}

In a case of reflexive bimodule a successive 
application of above canonical constructions give rise to the 
initial object \cite{AB}.

\section{Co-universal problem for Cartan pairs.}

Our aim in this section is to associate to any algebra $A$ its 
{\it co-universal} (right) Cartan pair.

Let $M,\, N$ be two $A$--bimodules and $\alpha: M\ra N$ a bimodule map 
between them.  Consider a transpose map $\alpha^T: N^*\ra M^*$ defined
by the formula
$$
<\alpha^T(Y),\, m>_M \doteq <Y,\, \alpha (m)>_N \eqno (3.1)
$$
where $Y\in N^*,\, m\in M$. From this definition one inspects that
$$
<\alpha^T(f.Y.g),\, m>_M = f\,<Y,\, \alpha (g. m)>_N =
f\,<\alpha^T(Y),\, m.g>_M = <f.\alpha^T(Y).g,\, m>_M 
$$
Thus we have proved

{\bf Proposition 3.1}. {\em
$\alpha^T$ is again an $A$-bimodule map.}

In what follows we shall apply above Proposition in order to dualise 
the universal differential calculus.

Recall, that  for a given algebra $A$ there exists the bimodule $\O$,
which is a kernel of the multiplication map and the differential 
calculus $(\O,\, d_u)$ with the following universal property:
for every differential calculus $(M,\, d)$ there exists a unique
bimodule map $\phi :\O\ra M$ such that the diagram\\
\ \ \ \ \ \\
$$
\matrix{
A &\mapright{d_u} &\O \cr
 &\stackrel{\ }{d\,}\,\searrow & \mapdown\phi \cr
 & & M \cr 
}
\eqno (3.2)
$$
\ \ \ \\
\ \ \ \ \\
is commutative. It is called the {\it universal differential calculus}
on $A$.

Put $\X\doteq (\O)^*$

{\bf Definition 3.2}. The right Cartan pair $(\X, \p_u)$ dual to
$(\O, d_u)$ is called a {\it right co-universal Cartan pair}.

Let us begin with the situation describe by the diagram (3.2).
Since, $\phi :\O\ra M$  is a bimodule map then its transpose 
$\phi^T : M^*\ra \X$ is again a bimodule map. Thus one has
$$
<X,\, df>_M = <X, \phi (d_u f)>_M = <\phi^T(X),\, d_u f>_{\O} \eqno (3.3)
$$
This means that $\phi^T(X)\in \X$ acts via $\p_u$ on $A$.

{\bf Theorem 3.3}. {\em  For any  algebra $A$ there exists
the (unique) co-universal right Cartan pair $(\X, \p_u)$. It has the 
following co-universal property: for an arbitrary right $A$-Cartan pair
$(N,\, \rho)$ there exists one and only one bimodule map $\Phi :\O\ra N $ 
such that the diagram\\
\ \ \ \\
$$
\matrix{
\X &\mapright{\p_u} &\nd_{\bk}(A) \cr
\mapup\Phi &\nearrow\,\stackrel{\,}{\p} &  \cr
N & &  \cr
}
\eqno (3.4)
$$
\  \ \ \\
\ \ \ \ \ \\
is commutative.}

The case when the Cartan pair $(N,\, \rho)$ originates from the calculus
$(M,\, d)$ has been considered above. The more general proof for an
arbitrary Cartan pair involves an explicit realization  
of $\X$. This will be done elsewhere.

\section{Algebra of differential operators.}

Let $(M, \p)$ be a (right) $A$-Cartan pair. To this data we are going to 
associate an algebra of (right handed, linear, partial) {\it differential 
operators} $\D (M, \p)$ together with its (algebraic) Fock space 
representation by means of operators acting on the algebra $A$.

More exactly, one has an algebra map
$l: A\ra \nd_{\bk}(A)$, which is an action of $A$ on $A$
$$
f^l (g) \doteq fg,\ \ \ \ \ (fg)^l = f^l\c g^l \eqno (4.1)$$
induced by the multiplication from the left. Intuitively, a differential 
operator is a polynomial expression in $X$ and $f$.
Thus, our task is to describe the subalgebra
$$
\D (M, \p)\doteq\hbox{gen}_{\nd_{\bk}}
\{f^l,\, \rp X\, | f\in A,\, X\in M \} \eqno (4.2)$$
generated in $\nd_{\bk}(A)$ by all endomorphism of the form $f^l, \rp X$.

{\bf Definition 4.1}. The algebra $\D (M, \p)$ is called an algebra
of a (linear, right handed) {\it differential operators} with respect 
to the right Cartan pair $(M, \p)$. 
Its elements are called {\it differential operators}.

We shall now try to understand better an algebraic structure of 
$\D (M, \p)$ in terms of an abstract algebra. Since $\p :M\ra \nd_{\bk}(A)$ 
is a $\bk$-linear map it uniquely lifts to an algebra map
$\hat\p :T_{\bk}M\ra\nd_{\bk}(A)$, where 
$T_{\bk}M\doteq\bk\oplus M\oplus (M\ot M)\ldots\equiv 
\bigoplus_{k\in\N}M^{\ot k} $ is a tensor algebra of 
the $\bk$-module $M$. 

Let $A\star T_{\bk}M$ denotes a {\it free product} of the algebras $A$
and $T_{\bk}M$ . This is the algebra generated by these algebras, with no
relation except for the identification of unite elements \cite{Avi}.
By the universal property of the free product $\star$\, , there exists 
a unique algebra map $\mu :A\star T_{\bk}M\ra \nd_{\bk} (A)$, such that 
it extends $l$ and $\hat \p$ {\it i.e.} 
$$
\mu|_A = l, \ \ \hbox{and}\ \ \ \mu|_{T_{\bk}M} = \hat\p \eqno (4.3)
$$
Thus by its very definition the algebra 
$$
\D (M, \p)\equiv A\star T_{\bk}M/\ker\mu \eqno (4.4)
$$
is isomorphic to the quotient algebra, where $\ker\mu$ denotes the
kernel of $\mu$. The Cartan pair axioms (2.4) and (2.5) rewritten 
in terms of the action read
$$
\rp{(f.X)} = f^l\c \rp X  \eqno (4.5)
$$
and
$$
\rp X\c f^l = \rp{(X.f)} + (\rp X (f))^l
\eqno (4.6)$$
This implies that $\ker\mu$ is nontrivial and contains an
ideal generated by the relations
$$
\{(f.X)-f\star X\,\ , (Y.g)-Y\star g +\rp Y (g)\,|\,f, g\in A; X, Y\in M \}
$$
In the case of modules over commutative algebra ({\it i.e.} $f.X=X.f$)
one sees that (4.5) - (4.6) reduce to the relations
$$
[\rp X\,,\, f^l] = (\rp X(f))^l \eqno (4.7)$$
which remains the classical {\it Canonical Commutation Relations} 
(CCR -- for short) between annihilation and creation operators. 
Here, $[a\,,\,b]\doteq a\c b - b\c a$ denotes as usual the 
commutator in $\nd_{\bk}(A)$.
In this sense the action of $\D (m, \p)$ on $A$ bears some features
of an (algebraic, oscillator type) Fock space representation: it extends 
the action $\p$ of $M$ on $A$ ({\it annihilation}) from one hand and the 
action of $A$ on $A$ ({\it creation}) from the other. The unit $1_A\in A$ 
plays the role of the vacuum (or ground) state ($\rp X(1_A)=0$ and 
$f^l(1_A)=f$). The commutation relations between creation operators are 
encoded in the algebra structure of $A$.
Observe that we have no {\it a priori} relations between the annihilation
operators: however, they may appear in the kernel of the action
$\hat \p :T_{\bk}M\ra \nd_{\bk}(A)$.

Obviously, the same construction can be applied to Cartan pairs
coming from differential calculi.
The true, ({\it i.e.} Hilbert) Fock space representation of a 
$q$-deformed differential calculus on the quantum plane has been 
discovered in \cite{PW}.

In the classical case  (of differential calculus on a manifold) 
the definition 4.1 gives rise to the definition of the classical 
(linear, partial) differential operators \cite{Woj}. For the case
of $\R^n$, these operators form so called a {\it Weyl algebra} 
\cite{Str}, which is generated by the CCR. 
In this way $\D (M, \p)$ becomes a generalization of the Weyl algebra.

\section{Connection versus covariant derivative.}

Our aim here is to verify  a possible application  of Cartan pairs
in a theory of connections. Provided with the notion of vector fields
one can try to generalize the Koszul's  approach: the noncommutative 
version of covariant derivatives.

Let $M$ be an $A$-bimodule and $E$ be a left $A$-module. 
$M\ot_A E$ is again a left $A$-module which can be contracted with
$M^*$. More exactly one has a $\bk$-linear map
$$\
\ll\,,\gg\, : M^*\ot_A (M\ot_A E)\ra E $$
given by the formula 
$$\ll X,\, m\ot_A \xi\gg\, \doteq \, <X, m>.\,\xi \eqno (5.1)$$
for $X\in M^*$, $m\in M$ and $\xi\in E$.

Following \cite{Co} a left connection on $E$ with respect to
the  differential calculus  $(M, d)$ is a linear map
$\nabla :E\ra M\ot_A E$ satisfying 
$$
\nabla (f.\,\xi) = f.\,\nabla \xi + df\ot_A\xi \eqno (5.2)
$$
(In the original Connes approach one uses the universal differential
calculus.)\\
Define 
$$
\nabla_X \xi\, \doteq\, \ll X\,, \nabla\xi\gg  \eqno (5.3) $$
as a {\it covariant derivative} along the vector field $X$ with 
respect to the connection $\nabla$. Thus, the
formula (5.3) defines the action 
$$M^*\ni X\mapsto \nabla_X\in \nd_{\bk}(E)$$
of $M$ on $E$. This action has the following properties
$$
\nabla_{f.X}\xi = f.\nabla_X\xi  
\eqno (5.4)$$
and
$$
\nabla_X (f.\,\xi) = \rp X (f).\,\xi + \nabla_{X.f}\xi
\eqno (5.5)$$
which generalize the axioms (2.4) - (2.5) of a right Cartan pair. 
From the other hand they may serve as an axiomatic definition of the 
covariant derivative and hence the connection $\nabla$.

  \end{document}